\documentstyle[twoside,12pt]{article}
\title{The expansion of the confluent hypergeometric function on the positive real axis}
\author{\sc R. B.\ Paris \\
{\em Division of Computing and Mathematics,}\\
{\em Abertay University, Dundee DD1 1HG, UK}}
\begin{document}
\def\f#1#2{\mbox{${\textstyle \frac{#1}{#2}}$}}
\def\dfrac#1#2{\displaystyle{\frac{#1}{#2}}}
\def\boldal{\mbox{\boldmath $\alpha$}}
{\newcommand{\Sgoth}{S\;\!\!\!\!\!/}
\newcommand{\bee}{\begin{equation}}
\newcommand{\ee}{\end{equation}}
\newcommand{\lam}{\lambda}
\newcommand{\ka}{\kappa}
\newcommand{\al}{\alpha}
\newcommand{\th}{\theta}
\newcommand{\fr}{\frac{1}{2}}
\newcommand{\fs}{\f{1}{2}}
\newcommand{\g}{\Gamma}
\newcommand{\br}{\biggr}
\newcommand{\bl}{\biggl}
\newcommand{\ra}{\rightarrow}
\newcommand{\mbint}{\frac{1}{2\pi i}\int_{c-\infty i}^{c+\infty i}}
\newcommand{\mbcint}{\frac{1}{2\pi i}\int_C}
\newcommand{\mboint}{\frac{1}{2\pi i}\int_{-\infty i}^{\infty i}}
\newcommand{\gtwid}{\raisebox{-.8ex}{\mbox{$\stackrel{\textstyle >}{\sim}$}}}
\newcommand{\ltwid}{\raisebox{-.8ex}{\mbox{$\stackrel{\textstyle <}{\sim}$}}}
\renewcommand{\topfraction}{0.9}
\renewcommand{\bottomfraction}{0.9}
\renewcommand{\textfraction}{0.05}
\newcommand{\mcol}{\multicolumn}
\date{}
\maketitle
\pagestyle{myheadings}
\markboth{\hfill \sc R. B.\ Paris  \hfill}
{\hfill \sc Algebraic asymptotic expansion \hfill}
\begin{abstract}
The asymptotic expansion of the Kummer function ${}_1F_1(a; b; z)$ is examined as $z\to+\infty$ on the Stokes line $\arg\,z=0$.  The correct form of the subdominant algebraic contribution is obtained for non-integer $a$. Numerical results demonstrating the accuracy of the expansion are given.
\vspace{0.4cm}

\vspace{0.3cm}

\end{abstract}
\vspace{0.3cm}

\noindent $\,$\hrulefill $\,$

\vspace{0.2cm}

\begin{center}
{\bf 1. \  Introduction}
\end{center}
\setcounter{section}{1}
\setcounter{equation}{0}
\renewcommand{\theequation}{\arabic{section}.\arabic{equation}}
The confluent hypergeometric function ${}_1F_1(a;b;z)$ (or first Kummer function
also denoted by $M(a, b, z)$) is defined for complex parameters $a$ and $b$ by
\[{}_1F_1(a; b; z)=\sum_{n=0}^\infty \frac{(a)_n z^n}{(b)_n n!}=1+\frac{a\,z}{b\,1!}+\frac{a(a+1) \,z^2}{b(b+1)\,2!}+\ldots\ ,\]
provided $b\neq 0, -1, -2, \ldots\,$, where $(a)_n=a(a+1) \ldots (a+n-1)=\g(a+n)/\g(a)$ is Pochhammer's symbol for the rising factorial. The series converges absolutely for all finite $z$ and reduces to a polynomial in $z$ when $a=0, -1, -2, \ldots\ $. We exclude this last case from our asymptotic considerations.
The function ${}_1F_1(a; b; z)$ is entire in $z$ and is consequently completely described in $-\pi<\arg\,z\leq\pi$.

The behaviour of ${}_1F_1(a; b; z)$ for large $z$ and fixed parameters is exponentially large in $\Re (z)>0$ and algebraic in character in $\Re (z)<0$. 
The well-known asymptotic expansion of ${}_1F_1(a; b; z)$ for $|z|\ra\infty$ is  given by \cite[p.~328]{DLMF}
\bee\label{e11}
\frac{\g(a)}{\g(b)}\,{}_1F_1(a;b;z)\sim E(z)+H(ze^{\mp\pi i}),
\ee
where the formal exponential and algebraic asymptotic series $E(z)$ and $H(z)$ are defined by
\bee\label{e12}
E(z):=z^{a-b} e^z \sum_{j=0}^\infty \frac{(b-a)_j(1-a)_j}{j! z^j},\quad
H(z):=\frac{z^{-a}\g(a)}{\g(b-a)} \sum_{k=0}^\infty \frac{(a)_k (1+a-b)_k}{k!\,z^k}.
\ee
In \cite[p.~328]{DLMF}, the sectors of validity of (\ref{e11}) are given by $-\fs\pi+\epsilon\leq\pm\arg\,z\leq\f{3}{2}\pi-\epsilon$, where $\epsilon>0$ denotes an arbitrarily small quantity. In the asymptotic theory of functions of hypergeometric type (see, for example, \cite[\S 2]{P17}, \cite[\S 2.3]{PK}) the expansion (\ref{e11}) is given in the sector $-\pi<\arg\,z\leq\pi$, with the upper or lower sign being chosen according as $\arg\,z>0$ or $\arg\,z<0$, respectively. In the sense of Poincar\'e there is no inconsistency between these two sets of sectorial validity.

The exponential expansion $E(z)$ is dominant as $|z|\ra\infty$ in $\Re (z)>0$ and becomes oscillatory on the anti-Stokes lines $\arg\,z=\pm\fs\pi$, where it is of comparable magnitude
to the algebraic expansion. In $\Re (z)<0$, the exponential expansion is subdominant with the behaviour of ${}_1F_1(a;b;z)$ then controlled by $H(ze^{\mp\pi i})$. The negative real axis $\arg\,z=\pm\pi$ is a Stokes line where $E(z)$ is {\it maximally} subdominant. When $H(ze^{\mp\pi i})$ is optimally truncated at, or near, its least term, the exponential expansion undergoes a smooth but rapid transition in the neighbourhood of  $\arg\,z=\pm\pi$; see \cite[p.~67]{DLMF} and \cite[Chapter 6]{PK}. It is clear, however, that (\ref{e11}) cannot account correctly  for the exponentially small expansion on $\arg\,z=\pi$, since it predicts the exponentially small behaviour $(|z|e^{\pi i})^{a-b} e^{-|z|}$ as $|z|\ra\infty$. 
When $a-b$ is non-integer with $a$ and $b$ real, this is a complex-valued contribution whereas ${}_1F_1(a; b; -|z|)$ is real.

The same argument applies on the other Stokes line $\arg\,z=0$, where $E(z)$ is {\it maximally} dominant. The subdominant algebraic expansion undergoes a Stokes phenomenon as the positive real axis is crossed. According to the second set of validity conditions of (\ref{e11}), just above this ray the multiplicative factor in front of the algebraic expansion is $(ze^{-\pi i})^{-a}$, whereas just below this ray the factor is $(ze^{\pi i})^{-a}$. There will be a smooth transition (at fixed $|z|$) between these two expressions. The first set of validity conditions mentioned above is more confusing, since it appears that on $\arg\,z=0$ one has the choice of either $(ze^{-\pi i})^{-a}$ or $(ze^{\pi i})^{-a}$ for the multiplicative factors. In the case of real parameters, this predicts a complex-valued contribution from the subdominant expansion, when in fact it clearly must be real-valued.

Although such subdominant terms are negligible in the Poincar\'e sense, their inclusion can significantly improve the numerical accuracy in computations; see, for example, \cite[p.~76]{Ob}.
The details of the expansion of ${}_1F_1(a;b;z)$ on the negative real axis have been discussed in \cite{P13}.
In this note we complete this discussion by considering the correct form of the subdominant algebraic expansion as $z\to+\infty$ on the positive real axis.

\vspace{0.6cm}

\begin{center}
{\bf 2. \ The expansion for ${}_1F_1(a; b; x)$ as $x\ra+\infty$}
\end{center}
\setcounter{section}{2}
\setcounter{equation}{0}
\renewcommand{\theequation}{\arabic{section}.\arabic{equation}}
In \cite[Theorem 1]{P13}, the expansion of the Kummer function ${}_1F_1(a;b;-x)$ for $x\to+\infty$ 
was established, which took into account the correct asymptotic behaviour of the exponentially small contribution on the negative real axis. If we replace the parameter $a$ by $b-a$ in \cite[(2.11)]{P13}, we have
\[\frac{1}{\g(b)}\,{}_1F_1(b-a;b;-x)-\frac{x^{a-b}}{\g(a)} \sum_{j=0}^{m_0-1}\frac{(b-a)_j(1-a)_j}{j! x^j} \hspace{4cm}\]
\bee\label{e21}
=\frac{x^{-a} e^{-x}}{\g(b-a)}\bl\{ \cos \pi a \sum_{j=0}^{M-1}(-)^jA_j\,x^{-j}+\frac{2\sin \pi a}{\sqrt{2\pi x}}\sum_{j=0}^{M-1}(-)^jB_j\,x^{-j}+O(x^{-M})\br\},
\ee
where the coefficients $A_j$ are given by
\bee\label{e14}
A_j=\frac{ (a)_j(1+a-b)_j}{j!}\qquad(j\geq 0)
\ee
and $M$ is a positive integer. Here, the dominant algebraic expansion on the left-hand side of (\ref{e21}), which has been optimally truncated with index $m_o$ given by
\bee\label{e22}
m_0=x+\Re (2a-b)+\alpha,\qquad |\alpha|<1,
\ee
has been subtracted off from ${}_1F_1(b-a;b;-x)$. 
The expansion on the right-hand side represents the exponentially small contribution, where the coefficients 
$B_j$, which are related to the $A_j$, are specified below.
If $a=n$, where $n$ is a positive integer, the sum on the left-hand side of (\ref{e21}) consists of $n$ terms and so cannot be optimally truncated. In this case, the expansion (\ref{e21}) holds with the upper limit of the sum on the left-hand side replaced by $n-1$ and with no contribution from the sum involving $B_j$ (since $\sin \pi a=0$); see \cite[Theorem 2]{P13}.

The coefficients $B_j$ are defined by
\bee\label{e23}
B_j=\sum_{k=0}^j (-2)^{k} (\fs)_k\,A_{j-k}\,G_{2k,j-k}.
\ee
The coefficients $G_{k,j}$ appear in the expansion of the so-called {\it terminant function} $T_\nu(z)$, which is defined as a multiple of the incomplete gamma function by
$T_\nu(z):=\xi(\nu)\,\g(1-\nu,z)$, $\xi(\nu)= e^{\pi i\nu}\g(\nu)/(2\pi i)$.
When 
\bee\label{e23a}
\nu=m_o-2a+b=x-\Im (2a-b)+\alpha,
\ee
so that $\nu\sim x$ as $x\to+\infty$, we have the expansion on the Stokes line $\arg\,z=\pi$ of this function  given by \cite[\S 5]{O}
\bee\label{e24}
T_{\nu-j}(xe^{\pi i})=\frac{1}{2}-\frac{i}{\sqrt{2\pi x}}\bl\{\sum_{k=0}^{M-1} (\fs)_{k} G_{2k,j} (\fs x)^{-k}+O(x^{-M})\br\},
\ee
for $j=0, 1, 2, \ldots$ and positive integer $M$.
The coefficients $G_{k,j}$ are computed from the expansion 
\[\frac{\tau^{\gamma_j-1}}{1-\tau}\,\frac{d\tau}{dw}=-\frac{1}{w}+\sum_{k=0}^\infty G_{k,j}w^k,\qquad \fs w^2=\tau-\log\,\tau-1,\]
where the parameter $\gamma_j$ is specified by
\bee\label{e25}
\gamma_j=\nu-x-j=\alpha-j-\Im (2a-b)\qquad(0\leq j\leq N-1)
\ee
by (\ref{e23a}).
The branch of $w(\tau)$ is chosen such that $w\sim \tau-1$ as $\tau\ra 1$ and 
so that upon reversion of the $w$-$\tau$ mapping
\[\tau=1+w+\f{1}{3}w^2+\f{1}{36}w^3-\f{1}{270}w^4+\f{1}{4320}w^5+ \cdots\ ,\]
it is found with the help of {\it Mathematica} that the first five even-order coefficients $G_{2k,j}\equiv 6^{-2k} {\hat G}_{2k,j}$ are\footnote{There was a misprint in the first term in ${\hat G}_{6,j}$ in \cite{P13}, which appeared as $-3226$ instead of $-3626$. This was pointed out by T. Pudlik \cite{TP}. The correct value was used in the numerical calculations described in \cite{P13}.}
\begin{eqnarray}
{\hat G}_{0,j}\!\!&=&\!\!\f{2}{3}-\gamma_j,\qquad {\hat G}_{2,j}=\f{1}{15}(46-225\gamma_j+270\gamma_j^2-90\gamma_j^3), \nonumber\\
{\hat G}_{4,j}\!\!&=&\!\!\f{1}{70}(230-3969\gamma_j+11340\gamma_j^2-11760\gamma_j^3+5040\gamma_j^4
-756\gamma_j^5),\nonumber\\
{\hat G}_{6,j}\!\!&=&\!\!\f{1}{350}(-3626-17781\gamma_j+183330\gamma_j^2-397530\gamma_j^3+370440\gamma_j^4
-170100\gamma_j^5\nonumber\\
&&\hspace{7cm}+37800\gamma_j^6-3240\gamma_j^7),\nonumber\\
{\hat G}_{8,j}\!\!&=&\!\!\f{1}{231000}(-4032746+43924815\gamma_j+88280280\gamma_j^2-743046480\gamma_j^3\nonumber\\
&&+1353607200\gamma_j^4-1160830440\gamma_j^5+541870560\gamma_j^6
-141134400\gamma_j^7\nonumber\\
&&\hspace{6cm}+19245600\gamma_j^8-1069200\gamma_j^9).\label{e26}
\end{eqnarray}
From this it is evident that the coefficients $B_j$ not only depend on $a$ and $b$ but also on $\alpha$ in (\ref{e23a}), which in turn depends on the particular value of the variable $x$ under consideration.

If we apply Kummer's transformation
\[{}_1F_1(a;b;z)=e^z\,{}_1F_1(b-a;b;-z)\]
to the hypergeometric function in (\ref{e21}) we obtain
\newtheorem{theorem}{Theorem}
\begin{theorem}$\!\!\!.$\ When $a$ is non-integer we have the expansion 
\[\frac{1}{\g(b)}\,{}_1F_1(a; b; x)-\frac{x^{a-b}e^x}{\g(a)} \sum_{j=0}^{m_o-1}\frac{(b-a)_j(1-a)_j}{j! \,x^j}\hspace{6cm}\]
\bee\label{e27}
=\frac{x^{-a}}{\g(b-a)} \bl\{ \cos \pi a \sum_{j=0}^{M-1}(-)^jA_j\,x^{-j}+\frac{2\sin \pi a}{\sqrt{2\pi x}}\sum_{j=0}^{M-1} (-)^jB_j\,x^{-j}+O(x^{-M})\br\}
\ee
as $x\ra+\infty$, where $m_o$ is the optimal truncation index  defined in (\ref{e22}) satisfying $m_o\sim x$ and $M$ is a positive integer. The coefficients $A_j$ and $B_j$ are given in (\ref{e14}) and (\ref{e23}).
\end{theorem}
The expansion on the right-hand side of (\ref{e27}) represents the correct form of the subdominant algebraic expansion
of ${}_1F_1(a;b;x)$ as $x\to+\infty$ on the Stokes line $\arg\,x=0$. 

When $a=n$, where $n$ is a positive integer, the sum on the left-hand side of (\ref{e27}) terminates after $n$ terms (and so cannot be optimally truncated) and the 
sum on the right-hand side involving the coefficients $B_j$ vanishes. In this case, the expansion (\ref{e27}) reduces to the standard result in (\ref{e11}), where there is no ambiguity caused by the leading factor $(e^{\mp\pi i}x)^{-a}$
appearing in front of the algebraic expansion.
\vspace{0.6cm}

\begin{center}
{\bf 3.\  Numerical examples and concluding remarks}
\end{center}
\setcounter{section}{3}
\setcounter{equation}{0}
\renewcommand{\theequation}{\arabic{section}.\arabic{equation}}
In this section we present some numerical examples to demonstrate the accuracy of the expansion in 
(\ref{e27}). As has already been noted, the coefficients $B_j$ depend on the parameters $a$ and $b$ and also on $\alpha$ (see the definition of $\gamma_j$ in (\ref{e25})), which appears in the value of the optimal truncation index $m_o$ in (\ref{e22}). The value of $\alpha$ clearly is a function of the particular value of $x$ being considered.
In Table 1 we show values of the coefficients $B_j$ for different $a$, $b$ and $\alpha$ (based on an integer value of $x$). In the calculations we have evaluated the coefficients $G_{2k,j}$ for $0\leq k\leq 6$.
\begin{table}[h]
\caption{\footnotesize{Values of the coefficients $B_j$ (with $\gamma_j=\alpha-j-\Im (2a-b)$) appearing in the expansion (\ref{e27}). The parameter $\alpha$ in (\ref{e22}) has been based on an integer value of $x$.}}
\begin{center}
\begin{tabular}{l|c|c|c}
\hline
&&&\\[-0.30cm]
%\mcol{1}{c|}{} & \mcol{2}{c|}{${}_1F_1(a;b;-x)$} & \mcol{1}{c}{$U(a,b,xe^{\pi i})$}\\
\mcol{1}{c|}{$j$} & \mcol{1}{c|}{$a=1/3,\ b=1$} & \mcol{1}{c|}{$a=3/4,\ b=1/2$} & \mcol{1}{c}{$a=1/4,\ b=3/4$}\\
\mcol{1}{c|}{} & \mcol{1}{c|}{$\alpha=1/3$} & \mcol{1}{c|}{$\alpha=0$} & \mcol{1}{c}{$\alpha=1/4$}\\
[.1cm]\hline
&&&\\[-0.30cm]
0 & 0.33333333333 & 0.66666666667 & 0.41666666667\\
1 & 0.15246913580 & 1.47731481481 & 0.16741898148\\
2 & 0.15301905742 & 3.83399470899 & 0.18232822317\\
3 & 0.35298391333 & 13.3434283401 & 0.44005533361\\
4 & 1.17133694166 & 59.8043967080 & 1.48497131914\\
5 & 4.97440127546 & 328.231235851 & 6.40730231996\\
6 & 26.0263928534 & 2130.38975509 & 33.9936918351\\
[.2cm]\hline
\end{tabular}
\end{center}
\end{table}

In Table 2 we show the values of the quantity ${\cal F}(x)$ defined by the left-hand side of (\ref{e27})
\[
{\cal F}(x):=\frac{1}{\g(b)}\,{}_1F_1(a; b; x)-\frac{x^{a-b}e^x}{\g(a)} \sum_{j=0}^{m_o-1} \frac{(b-a)_j (1-a)_j}{j!\,x^j}\]
compared with the truncated right-hand side
\[{\cal H}_M(x):=\frac{x^{-a}}{\g(b-a)} \bl\{ \cos \pi a \sum_{j=0}^{M}(-)^jA_j\,x^{-j}+\frac{2\sin \pi a}{\sqrt{2\pi x}}\sum_{j=0}^{M} (-)^jB_j\,x^{-j}\br\}\]
for different truncation index $M$.
In each case the optimal truncation index $m_o$ of the exponential expansion in ${\cal F}(x)$ was determined by inspection and the value of $\alpha$ determined from (\ref{e22}). In the case with $a=\fs$ 
the sum involving the coefficients $A_j$ makes no contribution to ${\cal H}_M(x)$. 
%When $a=2$, there is no optimal truncation index and Theorem 2 applies. 
It can be seen that the computed values of ${\cal F}(x)$ agree well with the  subdominant algebraic expansion.

\begin{table}[h]
\caption{\footnotesize{Values of ${\cal F}(x)$ and ${\cal H}_M(x)$ 
using an optimal truncation $m_o\sim x$ of the exponential expansion for different index $M$ in the subdominant algebraic expansion and parameters $a$, $b$ when $x=20$.}}
\begin{center}
\begin{tabular}{c||c|c|c}
\mcol{1}{c||}{} & \mcol{1}{c|}{$a=0.75,\ b=0.50$} & \mcol{1}{c|}{$a=0.50,\ b=1.25$} & \mcol{1}{c}{$a=-0.75,\ b=1.25$}\\
\mcol{1}{c||}{} & \mcol{1}{c|}{$\alpha=0,\ m_o=21$} & \mcol{1}{c|}{$\alpha=-0.75,\ m_o=19$} & \mcol{1}{c}{$\alpha=0.75,\ m_o=18$}\\
\mcol{1}{c||}{$M$} &\mcol{1}{c|}{${\cal H}_M(x)$} & \mcol{1}{c|}{${\cal H}_M(x)$} & \mcol{1}{c}{${\cal H}_M(x)$}\\
\hline
&&&\\[-0.25cm]
0 & 0.01343919981 & 0.04612049350 & $-6.58797674097$\\
1 & 0.01292521217 & 0.04684837963 & $-6.29520877265$\\
2 & 0.01296951096 & 0.04689194291 & $-6.29471792653$\\
3 & 0.01296356587 & 0.04689093247 & $-6.29471024690$\\
4 & 0.01296463732 & 0.04689125453 & $-6.29471013827$\\
5 & 0.01296439560 & 0.04689118984 & $-6.29471013739$\\
6 & 0.01296446100 & 0.04689120714 & $-6.29471013740$\\
&&&\\[-0.30cm]
\hline
&&&\\[-0.30cm]
${\cal F}(x)$ & 0.01296444571 & 0.04689120311  & $-6.29471013740$\\
\hline
\end{tabular}
\end{center}
\end{table}

We remark that the so-called Stokes multiplier on the positive real axis (given by the quantity in curly braces on the left-hand side of (\ref{e27})) is equal to $\cos \pi a$ to leading order. From (\ref{e11}) and (\ref{e12}) and the second set of validity conditions this quantity has the values $e^{\pi ia}$ and $e^{-\pi ia}$ just above and below $\arg\,z=0$. A commonly adopted, heuristic rule is that the Stokes multiplier on the Stokes line is given by the average of these two values to leading order, namely $\cos \pi a$. This agrees with the result stated in (\ref{e27}); see also \cite[p.~248]{PK}. 

It would be of interest to extend the result of Theorem 1 to the more general Wright function defined by
\[{}_1\Psi_1\bl(\!\!\begin{array}{c}(\alpha,a)\\(\beta,b)\end{array}\!\bl|z\br)=\sum_{n=0}^\infty \frac{\g(\alpha n+a)}{\g(\beta n+b)}\,\frac{z^n}{n!}\qquad (\kappa>0,\ |z|<\infty),\]
where $\alpha$, $\beta>0$ with $\kappa=1+\beta-\alpha$ and $a$, $b$ finite complex constants.
From the asymptotic theory of functions of this type (see \cite{P17}, \cite[\S\S 2.2.4, 2.3]{PK}) we find that
\[{}_1\Psi_1\bl(\!\!\begin{array}{c}(\alpha,a)\\(\beta,b)\end{array}\!\bl|z\br)\sim A Z^{\vartheta} e^Z \sum_{j=0}^\infty c_j Z^{-j}+\frac{1}{\alpha}\sum_{k=0}^\infty \frac{(-)^k\g(\frac{k+a}{\alpha})}{k!\g(b-\frac{\beta(k+a)}{\alpha})}\,(ze^{\mp \pi i})^{-(k+a)/\alpha}\]
as $|z|\to\infty$ in the sector $|\arg\,z|< \min\{\pi, \pi\kappa\}$, with the upper or lower sign chosen according as $\arg\,z>0$ or $\arg\,z<0$, respectively.
The various quantities in this asymptotic formula are given by
\[Z=\kappa(hz)^{1/\kappa},\quad h=\alpha^\alpha \beta^{-\beta},\quad \vartheta=a-b,\quad A=\kappa^{-\vartheta-\frac{1}{2}} \alpha^{a-\frac{1}{2}} \beta^{\frac{1}{2}-b}\]
and the $c_j$ are coefficients with $c_0=1$ and 
\[c_1=\frac{1}{12\alpha\beta}\bl\{\alpha(\alpha-1)(1-6b+6b^2)+\beta(\beta+1)(1-6a+6a^2)\]
\[+\alpha\beta[\alpha-\beta-2(1+6ab-6b)]\br\}.\]

Upon application of the above heuristic rule, we would expect the Stokes multiplier associated with the algebraic expansion on $\arg\,z=0$ to be equal to $\cos (\pi a/\alpha)$ to leading order. 
In the simpler case $\alpha=\beta$, a possible approach to investigate this conjecture is the integral representation
\[{}_1\Psi_1\bl(\!\!\begin{array}{c}(\alpha,a)\\(\alpha,b)\end{array}\!\bl|z\br)=\frac{1}{\g(b-a)} \int_0^1 t^{a-1}(1-t)^{b-a-1} \exp\,[zt^\alpha]\,dt\]
for $\Re (b)>\Re (a)>0$,
which reduces to that of $(\g(a)/\g(b))\,{}_1F_1(a;b;z)$ when $\alpha=1$.

\end{document}